\documentclass[10pt,a4paper,leqno,twoside]{article}
\usepackage{pifont}
\usepackage{bbding}
\usepackage{mathrsfs}
\usepackage{amsfonts,amssymb,amsmath,indentfirst,amsthm}
\usepackage{epsfig,cite}
\usepackage{color}

\numberwithin{equation}{section}
\newenvironment{abs}{\textbf{Abstract}\mbox{  }}{ }
\newenvironment{key words}{\emph{\texttt{Keywords}}\mbox{  }}{ }
\newtheorem{theorem}{Theorem}[section]
\newtheorem{lemma}[theorem]{Lemma}

\newtheorem{proposition}[theorem]{Proposition}

\renewenvironment{proof}{\noindent{\textbf{Proof.}}}{\hfill$\Box$}
\newenvironment{altproof}[1]
{\addvspace{0.3cm}\noindent
{\textbf{Proof of {#1}}.}}
{\nopagebreak\mbox{}\hfill $\Box$\par\addvspace{0.5cm}}
\theoremstyle{remark}
\newtheorem{remark}[theorem]{Remark}

\theoremstyle{plain}

\makeatletter
    
    \newcommand{\Rmnum}[1]{\expandafter\@slowromancap\romannumeral #1@}
\makeatother

\newcommand{\wt}{\widetilde}

\newcommand{\ov}{\overline}

   \def\dist{\mathrm{dist}}
   \def\span{\mathrm{span}}

   \def\N{\mathbb{N}}
   \def\R{\mathbb{R}}
   
   \def\Z{\mathbb{Z}}

\newcommand{\cO}{{\mathcal O}}

\newcommand{\al}{\alpha}
\newcommand{\be}{\beta}
\newcommand{\ga}{\gamma}
\newcommand{\de}{\delta}

\newcommand{\la}{\lambda}
\newcommand{\si}{\sigma}

\newcommand{\De}{\Delta}

\newcommand{\La}{\Lambda}
\newcommand{\Om}{\Omega}

\newcommand{\eps}{\varepsilon}

\newcommand{\pa}{\partial}

\def\Ker{\mathrm{Ker}}

\begin{document}

\title{\textbf{Multi-bubble nodal solutions to slightly subcritical elliptic problems with Hardy terms in symmetric domains}}
\author{Thomas Bartsch\thanks{Mathematisches Institut, Justus-Liebig-Universit$\ddot{a}$t Giessen, Arndtstr. 2, 35392 Giessen, Germany}, Qianqiao Guo\thanks{
 School of Mathematics and Statistics, Northwestern Polytechnical University, 710129 Xi'an, China}
}
\date{}
\maketitle

\noindent
\begin{abs}
We consider the slightly subcritical elliptic problem with Hardy term
\[
\left\{
\begin{aligned}
-\De u-\mu\frac{u}{|x|^2} &= |u|^{2^{\ast}-2-\eps}u &&\quad \text{in } \Om\subset\R^N, \\\
u &= 0&&\quad \text{on } \pa \Om,
\end{aligned}
\right.
\]
where $0\in\Om$ and $\Om$ is invariant under the subgroup $SO(2)\times\{\pm E_{N-2}\}\subset O(N)$; here $E_n$ denots the $n\times n$ identity matrix. If $\mu=\mu_0\eps^\al$ with $\mu_0>0$ fixed and $\al>\frac{N-4}{N-2}$ the existence of nodal solutions that blow up, as $\eps\to0^+$, positively at the origin and negatively at a different point in a general bounded domain has been proved in \cite{BarGuo-ANS}. Solutions with more than two blow-up points have not been found so far. In the present paper we obtain the existence of nodal solutions with a positive blow-up point at the origin and $k=2$ or $k=3$ negative blow-up points placed symmetrically in $\Om\cap(\R^2\times\{0\})$ around the origin provided a certain function $f_k:\R^+\times\R^+\times I\to\R$ has stable critical points; here $I=\{t>0:(t,0,\dots,0)\in\Om\}$. If $\Om=B(0,1)\subset\R^N$ is the unit ball centered at the origin we obtain two solutions for $k=2$ and $N\ge7$, or $k=3$ and $N$ large. The result is optimal in the sense that for $\Om=B(0,1)$ there cannot exist solutions with a positive blow-up point at the origin and four negative blow-up points placed on the vertices of a square centered at the origin. Surprisingly there do exist solutions on $\Om=B(0,1)$ with a positive blow-up point at the origin and four blow-up points on the vertices of a square with alternating positive and negative signs. The results of our paper show that the structure of the set of blow-up solutions of the above problem offers fascinating features and is not well understood.
\end{abs}

\vspace{.2cm}
\noindent
{\small\bf 2010 Mathematics Subject Classification:} 35B44, 35B33, 35J60.

\vspace{.2cm}
\noindent
{\small\bf Key words:}
Hardy term; Critical exponent; Slightly subcritical problems; Nodal solutions; Multi-bubble solutions.

\section{\textbf{Introduction}\label{Section 1}}
The paper is concerned with the semilinear singular problem
\begin{equation}\label{pro}
\left\{
\begin{aligned}
-\De u-\mu\frac{u}{|x|^2} &= |u|^{2^{\ast}-2-\eps}u &&\quad \text{in }\Om, \\\
u &= 0 &&\quad \text{on }\pa \Om,
\end{aligned}
\right.
\end{equation}
where $\Om\subset\R^N$, $N\ge7$, is a smooth bounded domain with $0\in\Om$; $2^*:=\frac{2N}{N-2}$ is the critical Sobolev exponent. In \cite{BarGuo-ANS}, we obtained the existence of two-bubble nodal solutions to problem \eqref{pro} that blow up positively at the origin and negatively at a different point in a general bounded domain, as $\eps\to0^+$ and $\mu=\mu_0\eps^\al$ with $\mu_0>0$ and $\al>\frac{N-4}{N-2}$. The location of the negative blow-up point is determined by the geometry of the domain. The existence of nodal bubble tower solutions has been proved in \cite{BarGuo-SNPDEA}. These are superpositions of positive and negative bubbles with different scalings, all blowing up at the origin. It seems to be a difficult and open problem whether solutions with a blow-up point at the origin and more than one blow-up point outside the origin exist in a general domain. In the present paper we investigate this problem in symmetric domains, in particular for the model case of the ball $\Om=B(0,1)$. 

In the case of a ball it is natural to place one blow-up point, say positive, at the origin and $k$ blow-up points, say negative, at the vertices of a regular $k$-gon with center at the origin. We shall prove that solutions of this shape exist for $2\le k\le3$ but, somewhat surprisingly, not for $k=4$. On the other hand, we prove the existence of solutions with four blow-up points, two positive and two negative ones, at the vertices of a square, centered symmetrically around the positive blow-up point at the origin. Our results show that the existence of solutions of \eqref{pro} with three or more blow-up points is interesting and far from being understood.

When $\mu=0$ the blow-up phenomenon for positive and for nodal solutions to problem \eqref{pro} has been studied extensively, see for instance \cite{BaLiRey95CVPDE,BarDPis-AIHPNA,BarDPis-BLMS,BarMiPis06CVPDE, BrePe89PNDEA,PiDolMusso03JDE,FluWeip7MM,GroTaka10JFA,Han91AIHPAN, MussoPis10JMPA,PisWe07AIHPA,Rey89MM,Rey90JFA,Rey91DIE} and the references therein. However for $\mu\neq0$, few results are known about the existence of positive or nodal solutions with multiple bubbles to problem \eqref{pro}. Positive solutions have been obtained in \cite{Cao-Peng:2006}. Related results, though for different equations, can be found in \cite{FelliPis06CPDE,FelliTer05CCM,MussoWei12IMRN}. We also want to mention the papers \cite{CaoH04JDE,CaoP03JDE,EkeGhou02,FG,GY,GuoNiu08JDE,Jan99,RuWill03,Smets05TAMS,Terr96} dealing with the critical exponent, i.e.\ $\eps=0$.


%

An important role will be played by the limit problem
\begin{equation}\label{eq:limit1}
\left\{
\begin{aligned}
-\De u-\mu\frac{u}{|x|^{2}} &= |u|^{2^{\ast}-2}u &&\quad \text{in } \R^N, \\
u &\to 0 &&\quad \text{as } |x|\to\infty
\end{aligned}
\right.
\end{equation}
which has been investigated in \cite{CatW00,Terr96}. Positive solutions of \eqref{eq:limit1} in the range $0\le\mu<\ov{\mu}:=\frac{(N-2)^2}{4}$ are given by
\[
  V_{\mu,\si}(x) = C_\mu\left(\frac{\si}{\si^2|x|^{\be_1}+|x|^{\be_2}}\right)^{\frac{N-2}{2}}
\]
with $\si>0$, $\be_1:=(\sqrt{\ov{\mu}}-\sqrt{\ov{\mu}-\mu})/\sqrt{\ov{\mu}}$, $\be_2:=(\sqrt{\ov{\mu}}+\sqrt{\ov{\mu}-\mu})/\sqrt{\ov{\mu}}$, and
$C_\mu:=\left(\frac{4N(\ov{\mu}-\mu)}{N-2}\right)^{\frac{N-2}{4}}$.
These solutions minimize
\[
  S_\mu
   := \min_{u\in D^{1,2}(\R^N)\setminus \{0\}}
       \frac{\int_{\R^N}(|\nabla u|^2-\mu\frac{u^2}{|x|^2})dx}
            {(\int_{\R^N}|u|^{2^*}dx)^{2/{2^*}}},
\]
and there holds
$$
  \int_{\R^N} \left(|\nabla V_{\mu,\si}|^2-\mu\frac{|V_{\mu,\si}|^2}{|x|^2}\right)dx
    = \int_{\R^N} |V_{\mu,\si}|^{2^*}dx = S_\mu^{\frac{N}{2}}.
$$
In the range $0<\mu<\ov{\mu}$ these are all positive solutions of \eqref{eq:limit1}. In the case $\mu=0$ these are all solutions with maximum at $x=0$. Of course, if $\mu=0$ also translates of $V_{\mu,\si}$ are solutions of
\begin{equation}\label{eq:limit2}
\left\{
\begin{aligned}
  -\De u &= |u|^{2^{\ast}-2}u &&\quad \text{in } \R^N, \\\
  u &\to 0 &&\quad \text{as } |x|\to\infty.
\end{aligned}
\right.
\end{equation}
We will write
\[
  U_{\de,\xi}(x) = C_0\left(\frac{\de}{\de^2+|x-\xi|^2}\right)^{\frac{N-2}{2}}
\]
for the solutions of \eqref{eq:limit2} where $\de>0$, $\xi\in \R^N$ and $C_0:=(N(N-2))^{\frac{N-2}{4}}$. These are the well known Aubin-Talenti instantons (see \cite{Aub76,Tal76}).

Now we state our main results. We consider domains satisfying the condition
\begin{itemize}\label{eq:domain}
\item[$(A_1)$] $\Om\subset\R^N$ is a bounded domain with $0\in\Om$, and it is invariant under the subgroup $SO(2)\times\{\pm E_{N-2}\}\subset O(N)$.
\end{itemize}
We use the notation $x=(x',x'')\in\Om\subset\R^2\times\R^{N-2}$ and write $A(x',x'')=(Ax',x'')$ for $A\in SO(2)$. For $k\in\N$ let $R_k=\begin{pmatrix}\cos\frac{2\pi}{k} & -\sin\frac{2\pi}{k} \\ \sin\frac{2\pi}{k} & \cos\frac{2\pi}{k} \end{pmatrix} \in SO(2)$, and set $I=\{t>0:(t,0,\dots,0)\in\Om\}\subset\R$. 

Our first results are concerned with the existence of nodal solutions with $k+1$ bubbles, one being positive and $k$ being negative. Let $G(x,y)=\frac{1}{|x-y|^{N-2}}-H(x,y)$ be the Green function (up to a coefficient involving the volume of the unit ball) for the Dirichlet Laplace operator in $\Om$ with regular part $H$. For $k=2,3$ we define the function $f_k:\R^+\times\R^+\times I \to \R$ by
\[
\begin{aligned}
f_k(\la_0,\la_1,t)
 &:= b_1\bigg(H(0,0)\la_0^{N-2} + kH(\xi(t),\xi(t))\la_1^{N-2}
        + 2kG(\xi(t),0)\la_0^{\frac{N-2}{2}}\la_1^{\frac{N-2}{2}}\\
 &\hspace{1cm}   - 2\binom{k}{2}G(\xi(t),R_k\xi(t))\la_1^{N-2}\bigg)- b_2\frac{N-2}{2}\ln\left(\la_0\la_1^k\right).
\end{aligned}
\]
where $\xi(t)=(t,0,\dots,0)$ and
\[
  b_1 = \frac{1}{2}C_0\int_{\R^N}U_{1,0}^{2^*-1}\quad\text{and}\quad b_2=\frac{1}{2^*}\int_{\R^N} U_{1,0}^{2^*}.
\]
Finally we call a critical point of $f_k$ stable if it is isolated and has nontrivial local degree. This is the case, for instance, if it is non-degenerate or an isolated local maximum or minimum.

\begin{theorem}\label{thm:k bubbles}
Suppose $\Om$ satisfies $(A_1)$, and suppose $(\la_0,\la_1,t)\in\R^+\times\R^+\times I$ is a stable critical point of $f_k$, $k=2$ or $k=3$. Let $\mu_0>0$ and $\al>\frac{N-4}{N-2}$ be fixed. Then there exists $\eps_0>0$ such that for every $\eps\in(0,\eps_0)$ problem \eqref{pro}  with $\mu=\mu_\eps=\mu_0\eps^\al$ has a pair of solutions $\pm u_\eps$ satisfying
\begin{equation}\label{ball-shape of solution-k bubbles}
\begin{aligned}
u_\eps(x)
 &= V_{\mu_\eps,\si_\eps}(x) - \sum_{i=1}^k U_{\de_{\eps},R_k^i\xi(t_\eps)}(x) + o(1)\\
 &= C_{\mu_\eps}\left(\frac{\si_\eps}{(\si_\eps)^2|x|^{\be_1}
       + |x|^{\be_2}}\right)^{\frac{N-2}{2}}
    - C_0\sum_{i=1}^k\left(\frac{\de_\eps}{(\de_\eps)^2
       + |x-R_k^i\xi(t_\eps)|^2}\right)^{\frac{N-2}{2}}
    + o(1),
\end{aligned}
\end{equation}
where 
\[
  \si_\eps = \big(\la_0+o(1)\big)\eps^{\frac{1}{N-2}},\ \de_\eps = \big(\la_1+o(1)\big)\eps^{\frac{1}{N-2}},\ 
  \xi(t_\eps) = (t_\eps,0,\dots,0)= (t+o(1),0,\dots,0)\quad\text{as $\eps\to0$.}
\]
These solutions satisfy the following symmetries:
\begin{equation}\label{eq:symmetry of solutions 1}
  u_\eps(x',x'') = u_\eps(x',-x'') = u_\eps(R_kx',x'') \quad\text{for $(x',x'')\in\Om\subset\R^2\times\R^{N-2}$}.
\end{equation}
\end{theorem}

As Proposition~\ref{prop:k=4} below shows, for $k=4$ a family $u_\eps$ as in Theorem~\ref{thm:k bubbles} need not exist even in the case of the ball. It is a challenging problem to find critical points of $f_k$ for general domains satisfying $(A_1)$. We consider the special case where $\Om=B(0,1) \subset \R^N$ is the unit ball.

\begin{theorem}\label{theorem-ball-k bubbles}
If $\Om=B(0,1)\subset\R^N$, $k=2$ and $N\ge7$, or $k=3$ and $N$ is large enough, then $f_k$ has two stable critical points, one is a local minimum, the other a mountain pass point with Morse index 1. As a consequence, problem \eqref{pro} has two families of solutions $\pm u^1_\eps$, $\pm u^2_\eps$ as in Theorem~\ref{thm:k bubbles}. They have the additional symmetry
\[
  u_\eps(x',x'') = u_\eps(x',Ax'')\quad\text{for all $A\in SO(N-2)$.}
\]
\end{theorem}

\begin{remark}
a) In the proof of the case $k=3$ we provide an explicit inequality for $N$, so that the solutions as in Theorem~\ref{theorem-ball-k bubbles} exist if this inequality holds. Numerical computations show that this inequality is not satisfied for $N=7$. We do not know the optimal value for $N$ such that Theorem~\ref{theorem-ball-k bubbles} is true for $k=3$; see also Remark~\ref{rem:nonexist}.

b) We conjecture that Theorem~\ref{theorem-ball-k bubbles} holds for other domains satisying $(A_1)$, for instance for $\Om=B^2(0,1)\times\Om' \subset \R^2\times\R^{N-2}$ with $\Om' = -\Om'\subset\R^{N-2}$ a bounded symmetric neighborhood of $0$. Our proof of Theorem~\ref{theorem-ball-k bubbles} uses the explicit knowledge of the Green function for the ball, hence it does not extend immmediately to other domains. 
\end{remark}

The next result shows that Theorem~\ref{theorem-ball-k bubbles} is optimal.

\begin{proposition}\label{prop:k=4}
For $\Om=B(0,1)\subset\R^N$, $N\ge7$ and $k=4$ there does not exist a family of solutions $\pm u_\eps$ as in Theorem~\ref{thm:k bubbles}. 
\end{proposition}

\begin{remark}
a) We conjecture that Proposition~\ref{prop:k=4} can be generalized to $k\ge4$. 

b) We cannot exclude the existence of solutions with a positive bubble at the origin and $k=4$ negative bubbles placed somewhere in the ball and with possibly different blow-up parameters $\de$. However, we can show that there do not exist solutions with four negative bubbles at the vertices of a square centered at the origin even if one allows different blow-up speeds, i.e. if one replaces the $\de_\eps$ in \eqref{ball-shape of solution-k bubbles} by $\de_{i,\eps}$, $i=1,\dots,4$. In fact, it is not difficult to prove that the blow-up parameters $\de_{i,\eps}$ have to be independent of $i$ if the vertices are at $R_4^i\xi(t_\eps)$, $i=1,\dots,4$. 
\end{remark}

Considering Proposition~\ref{prop:k=4} the following existence results of nodal solutions with five bubbles, three being positive and two being negative, is somewhat surprising.

\begin{theorem}\label{theorem-ball-k bubbles-add}
Let $\Om=B(0,1)\subset\R^N$, $N\ge7$, $\mu_0>0$ and $\al>\frac{N-4}{N-2}$ be fixed. Then there exists $\eps_0>0$ such that for any $\eps\in(0,\eps_0)$, there exist a pair of $5$-bubble solutions $\pm u_\eps$ to problem \eqref{pro} with $\mu=\mu_\eps=\mu_0\eps^\al$ of the shape
\[
\begin{aligned}
u_\eps(x)
 &= V_{\mu_\eps,\si_\eps}(x) + \sum_{i=1}^4 (-1)^iU_{\de_{i,\eps},R_4^i\xi(t_\eps)}(x) + o(1)\\
 &= C_{\mu_\eps}\left(\frac{\si_\eps}{(\si_\eps)^2|x|^{\be_1} + |x|^{\be_2}}\right)^{\frac{N-2}{2}}
       + C_0\sum_{i=1}^4(-1)^i\left(\frac{\de_{i,\eps}}{(\de_{i,\eps})^2 + |x-R_4^i\xi(t_\eps)|^2}\right)^{\frac{N-2}{2}}
       + o(1),
\end{aligned}
\]
where $\si_\eps = \big(\la_0+o(1)\big)\eps^{\frac{1}{N-2}}$, $\de_{1,\eps} = \de_{3,\eps} = \big(\la_{1}+o(1)\big)\eps^{\frac{1}{N-2}}$, $\de_{2,\eps} = \de_{4,\eps} = \big(\la_{2}+o(1)\big)\eps^{\frac{1}{N-2}}$, $\xi(t_\eps) = (t_\eps,0,\dots,0)=(t+o(1),0\dots,0)$ as $\eps\to0$,
for some $\la_0,\la_1,\la_2>0$, $t\in (0,1)$. These solutions satisfy the symmetries:
\begin{equation}\label{eq:symmetry of solutions 2}
  u_\eps(x',x'') = u_\eps(x',Ax'') = -u_\eps(R_4x',x'') \quad\text{for $(x',x'')\in B(0,1)\subset\R^2\times\R^{N-2}$, $A\in SO(N-2)$}.
\end{equation}
\end{theorem}

\begin{remark}
a) The parameters $(\la_0,\la_1,\la_2,t)\in\R^+\times\R^+\times\R^+\times(0,1)$ in Theorem~\ref{theorem-ball-k bubbles-add} are obtained as a critical point of a suitable limit function $f_5$. We conjecture that there exists a second solution in Theorem~\ref{theorem-ball-k bubbles-add} but the computations for finding a second critical point of $f_5$ are intimidating.

b) It seems that for $k=2$ in Theorem \ref{theorem-ball-k bubbles}, it is still possible to obtain the information on the nodal sets of the solutions as in 
\cite{BarDPis-BLMS}. For $k=3$ in Theorem~\ref{theorem-ball-k bubbles} and for Theorem ~\ref{theorem-ball-k bubbles-add}, it is an interesting problem to study the profile of the nodal sets of the solutions.

c) As stated in \cite{BarGuo-ANS}, the assumption $\al>\frac{N-4}{N-2}$ is essential in our theorems.
\end{remark}

The paper is organized as follows. In Section~2, we collect some notations and preliminary results. Section~3 is devoted to the method of finite dimensional reduction. Section~\ref{sec:proof k=2,3} contains the proof of Theorems~\ref{thm:k bubbles} and \ref{theorem-ball-k bubbles}. Proposition~\ref{prop:k=4} is proved in Section~\ref{sec:proof prop}, and finally Theorem~\ref{theorem-ball-k bubbles-add} is proved in Section~\ref{sec:proof 5 bubbles}.

\section{Notations and preliminary results}\label{Section 2}

Throughout this paper, positive constants will be denoted by $C, c$. By Hardy's inequality the norm
\[
  \|u\|_\mu:=\left(\int_{\Om}(|\nabla u|^2-\mu\frac{u^2}{|x|^2})dx\right)^{\frac12}
\]
is equivalent to the norm $\|u\|_0=\left(\int_{\Om}|\nabla u|^2dx\right)^{1/2}$ on $H_0^1(\Om)$ provided $0\le\mu<\ov{\mu}$. This inequality is of course satisfied for $\mu=\mu_0\eps^\al$ with $\al>0$ and $\eps>0$ small. We write $H_\mu(\Om)$ for the Hilbert space consisting of the $H^1_0(\Omega)$ functions with the inner product
\[
  (u,v)_\mu:=\int_{\Om}\left(\nabla u\nabla v-\mu\frac{uv}{|x|^2}\right)dx.
\]
As in \cite{FelliPis06CPDE, BarGuo-ANS} let $\iota_\mu^*:L^{2N/(N+2)}(\Om)\to H_\mu(\Om)$ be the adjoint operator of the inclusion $\iota_\mu:H_\mu(\Om)\to L^{2N/(N-2)}(\Om)$, that is,
\begin{equation*}
\iota_\mu^*(u) = v\qquad\Longleftrightarrow\qquad
 (v,\phi)_\mu = \int_{\Om}u(x)\phi(x)dx,\quad\text{for all }\phi\in H_\mu(\Om).
\end{equation*}
There exists $c>0$ such that
\begin{equation*}
\|\iota_\mu^*(u)\|_{\mu}\le c\|u\|_{2N/(N+2)}.
\end{equation*}
Then problem \eqref{pro} is equivalent to the fixed point problem
\[
u=\iota_\mu^*(f_\eps(u)),\quad u\in H_\mu(\Om),
\]
where $f_\eps(s)=|s|^{2^*-2-\eps}s$.

The following proposition is from \cite[Proposition~3.1]{BarGuo-ANS}. 
\begin{proposition}\label{proposition-eigenvalue}
Let $0<\mu<\overline{\mu}$, and let $\La_i$, $i=1,2,\dots$, be the eigenvalues of
\begin{equation*}
\begin{cases}
-\De u-\mu\frac{u}{|x|^2} = \La|V_\si|^{2^{\ast}-2}u &\quad\text{in } \R^N,\\
|u|\to0 &\quad\text{as }|x|\to+\infty
\end{cases}
\end{equation*}
in increasing order. Then $\La_1=1$ with eigenfunction $V_\si$, $\La_2=2^*-1$ with eigenfunction $\frac{\pa V_\si}{\pa \si}$.
\end{proposition}

Our main results will be proved using variational and singular limit methods applied to the energy functional
\[
  J_\eps(u)
   := \frac12\int_{\Om}\left(|\nabla u|^2-\mu\frac{u^2}{|x|^2}\right)dx
     - \frac{1}{2^*-\eps}\int_{\Om}|u|^{2^*-\eps}dx
\]
defined on $H_\mu(\Om)$.

Let us also recall that the Green's function of the Dirichlet Laplacian $G(x,y)=\frac{1}{|x-y|^{N-2}}-H(x,y)$ and its regular part $H$ are symmetric: $G(x,y)=G(y,x)$ and $H(x,y)=H(y,x)$. If $\Om$ is invariant under some $A\in O(N)$ then $G(Ax,Ay)=G(x,y)$, and the same holds for $H$.

\section{The finite dimensional reduction}\label{section 3}

First we recall some notation from \cite{BarGuo-ANS}. Fix $\mu_0>0$, $\al>\frac{N-4}{N-2}$ and an integer $k\ge0$. For $\la=(\la_0,\la_1,\dots,\la_{k})\in \R_+^{k+1}$ and $\xi=(\xi_1,\xi_2,\dots,\xi_{k})\in\Om^{k}$ we define
\[
W_{\eps,\la,\xi}
 :=\sum_{i=1}^{k}\Ker\left(-\De-(2^*-1)U_{\de_i,\xi_i}^{2^*-2}\right)
     + \Ker\left(-\De-\frac{\mu_\eps}{|x|^2}-(2^*-1)V_{\mu_\eps,\si_\eps}^{2^*-2}\right) 
 \subset H^1(\R^N)
\]
where $\de_i=\la_i\eps^{\frac{1}{N-2}}$, $\mu_\eps=\mu_0\eps^{\al}$, $\si_\eps=\la_0\eps^{\frac{1}{N-2}}$. By Proposition \ref{proposition-eigenvalue} and \cite{BianEg91JFA} we know that
\[
W_{\eps,\la,\xi}
 = \span\left\{\Psi_i^j,\ \Psi_i^0,\ \Psi_0,\ i=1,2,\dots,k,\ j=1,2,\dots,N\right\},
\]
where for $i=1,2,\dots,k$ and $j=1,2,\dots,N$:
$$
\Psi_i^j:=\frac{\pa U_{\de_i,\xi_i}}{\pa \xi_{i,j}},\quad
\Psi_i^0:=\frac{\pa U_{\de_i,\xi_i}}{\pa \de_i},\quad
\Psi_0:=\frac{\pa V_{\mu_\eps,\si_\eps}}{\pa \si}
$$
with $\xi_{i,j}$ the $j$-th component of $\xi_i$. For $\eta\in(0,1)$ we define
\[
\begin{aligned}
\cO_\eta
 &:=\big\{(\la,\xi)\in\R_+^{k+1}\times \Om^{k}: \la_i\in(\eta,\eta^{-1}) \text{ for $i=0,\dots,k$},\ \dist(\xi_i,\pa\Om)>\eta,\\
 &\hspace{2cm}
 |\xi_i|>\eta,\ |\xi_{i_1}-\xi_{i_2}|>\eta,\ \text{for $i,i_1,i_2=1,\dots,k,\ i_1\neq i_2$}\big\}.
\end{aligned}
\]

The projection $P:H^1(\R^N)\to H_0^1(\Om)$ is defined by $\De Pu=\De u$ in $\Om$ and $Pu=0$ on $\pa\Om$. We also need the spaces
\[
K_{\eps,\la,\xi}:=P W_{\eps,\la,\xi}
\]
and
\[
K^\bot_{\eps,\la,\xi}
 :=\{\phi\in H_\mu(\Om):(\phi,P\Psi)_{\mu_\eps}=0,\text{ for all }\Psi\in W_{\eps,\la,\xi}\},
\]
as well as the $(\cdot,\cdot)_{\mu_\eps}$-orthogonal projections
\[
\Pi_{\eps,\la,\xi} :H_{\mu_\eps}(\Om)\to K_{\eps,\la,\xi}
\]
and
\[
\Pi^\bot_{\eps,\la,\xi} := Id-\Pi_{\eps,\la,\xi}:H_{\mu_\eps}(\Om)\to K^\bot_{\eps,\la,\xi}.
\]

Then solving problem \eqref{pro} is equivalent to finding $\eta>0$, $\eps>0$, $(\la,\xi)\in\cO_\eta$ and $\phi_{\eps,\la,\xi}\in K^\bot_{\eps,\la,\xi}$ such that:
\begin{equation}\label{main-equality Com-1}
\Pi^\bot_{\eps,\la,\xi}\left(V_{\eps,\la,\xi} + \phi_{\eps,\la,\xi}
 - \iota_\mu^*(f_\eps(V_{\eps,\la,\xi} + \phi_{\eps,\la,\xi}))\right)
= 0,
\end{equation}
and
\begin{equation*} 
\Pi_{\eps,\la,\xi}\left(V_{\eps,\la,\xi} + \phi_{\eps,\la,\xi} - \iota_\mu^*(f_\eps(V_{\eps,\la,\xi}
 + \phi_{\eps,\la,\xi}))\right)
= 0,
\end{equation*}
where in the case of Theorem~\ref{theorem-ball-k bubbles}
\begin{equation}\label{shape of V-1}
V_{\eps,\la,\xi}=-\sum\limits_{i=1}^{k}PU_{\de_i,\xi_i}+PV_{\mu_\eps,\si_\eps}
\end{equation}
with $k=2,3$, and in the case of Theorem~\ref{theorem-ball-k bubbles-add}
\begin{equation}\label{shape of V-add1}
V_{\eps,\la,\xi}=\sum\limits_{i=1}^{k}(-1)^iPU_{\de_i,\xi_i}+PV_{\mu_\eps,\si_\eps}
\end{equation}
with $k=4$.



The following two propositions have been proved in \cite{BarGuo-ANS}.

\begin{proposition}\label{proposition-estimate of error-1}
For every $\eta>0$ there exist $\eps_0>0$ and $c_0>0$ with the following property. For every $(\la,\xi)\in\cO_\eta$ and for every $\eps\in(0,\eps_0)$ there exists a unique solution $\phi_{\eps,\la,\xi}\in K^\bot_{\eps,\la,\xi}$ of equation \eqref{main-equality Com-1} satisfying
\begin{equation*}
\|\phi_{\eps,\la,\xi}\|_{\mu_\eps}
 \le c_0\left(\eps^{\frac{N+2}{2(N-2)}}+\eps^{\frac{1+2\al}{4}}\right).
\end{equation*}
The map $\Phi_\eps:\cO_\eta\to K^\bot_{\eps,\la,\xi}$ defined by $\Phi_\eps(\la,\xi):=\phi_{\eps,\la,\xi}$ is $C^1$.
\end{proposition}

Now we can define the reduced functional
\[
  I_\eps:\cO_\eta \to \R,\quad I_\eps(\la,\xi) := J_\eps(V_{\eps,\la,\xi}+\phi_{\eps,\la,\xi}).
\]

\begin{proposition}\label{proposition-reducement-1}
If $(\la,\xi)\in\cO_\eta$ is a critical point of $I_\eps$ then $V_{\eps,\la,\xi}+\phi_{\eps,\la,\xi}$ is a solution to problem \eqref{pro} for $\eps>0$ small.
\end{proposition}

So far everything works on a general bounded domain. Now we will use the invariance of $I_\eps$ under certain symmetries for further reductions. For $A\in O(N)$, $\xi=(\xi_1,\dots,\xi_k)\in(\R^N)^k$ and $u\in L^p(\R^N)$ we set $A\xi:=(A\xi_1,\dots,A\xi_k)$ and $A*u:=u\circ A^{-1}$. This induces isometric actions of $O(N)$ on $(\R^N)^k$ as well as on $L^p(\R^N)$ and, if $A\Om=\Om$,  on $L^p(\Om)$ and on $H_\mu(\Om)$ such that $\iota_\mu$ and $\iota_\mu^*$ are equivariant. Moreover we have
\[
  U_{\de,A\xi} = A*U_{\de,\xi}\quad\text{and}\quad W_{\eps,\la,A\xi} = \{A*u:u\in W_{\eps,\la,\xi}\},
\]
and analogously for $K_{\eps,\la,\xi}$, $\Pi_{\eps,\la,\xi}$, $V_{\eps,\la,\xi}$. As a consequence, the uniqueness statement in Proposition~\ref{proposition-estimate of error-1} implies
\begin{equation}\label{eq:sym phi}
  \phi_{\eps,\la,A\xi} = A*\phi_{\eps,\la,\xi},
\end{equation}
hence $I_\eps$ is invariant with respect to the action $A*(\la,\xi)=(\la,A\xi)$ of $O(N)$ on $\cO_\eta$:
\[
  I_\eps(\la,A\xi) = I_\eps(\la,\xi).
\]
Now we apply the principle of symmetric criticality using the matrix $A_N:=\begin{pmatrix}E_2 & 0 \\ 0 & -E_{N-2}\end{pmatrix}\in O(N)$. By assumption $A_N(\Om)=\Om$, hence $A_N$ acts on $\cO_\eta$ as above leaving $I_\eps$ invariant. The principle of symmetric criticality implies that critical points of $I_\eps$ constrained to the fixed point set
\[
  \cO_\eta^{A_N} = \{(\la,\xi)\in\cO_\eta: A_N\xi=\xi\} = \{(\la,\xi)\in\cO_\eta: \xi_i=(\xi_i',0)\in\R^2\times\R^{N-2},\ i=1,\dots,k\}
\]
are critical points of $I_\eps$. We also need the invariance of $I_\eps$ with respect to permutations of the blow-up points. Here we need to distinguish between the cases where $V_{\eps,\la,\xi}$ is of the form \eqref{shape of V-1} or of the form  \eqref{shape of V-add1}. Let $S_k$ denote the group of permutations of $\{1,\dots,k\}$. For $\pi\in S_k$ and $(\la,\xi)\in\R^{k+1}\times(\R^N)^k$ we define
\[
  \pi*(\la_0,\la_1,\dots,\la_k) := (\la_0,\la_{\pi(1)},\dots,\la_{\pi(k)})\quad\text{and}\quad \pi*(\xi_1,\dots,\xi_k) := (\xi_{\pi(1)},\dots,\xi_{\pi(k)}).
\]
In the case when $V_{\eps,\la,\xi}$ is of the form \eqref{shape of V-1} it is obvious that
\[
  I_\eps(\pi*\la,\pi*\xi) = I_\eps(\la,\xi)\quad\text{for all $\pi\in S_k$, $(\la,\xi)\in \cO_\eta$.}
\]
It follows that $I_\eps$ is invariant under the map
\[
  \tau:\cO_\eta^{A_N} \to \cO_\eta^{A_N},\quad 
  \tau(\la_0,\la_1,\dots,\la_k,\xi_1,\dots,\xi_k) := (\la_0,\la_k,\la_1,\dots,\la_{k-1},R_k\xi_k,R_k\xi_1,\dots,R_k\xi_{k-1}),
\]
which induces an action of $\Z/k\Z$ on $\cO_\eta^{A_N}$; here $R_k(\xi',\xi''):=(R_k\xi',\xi'')$ where $R_k\in SO(2)$ is the rotation from Theorem~\ref{theorem-ball-k bubbles}. Therefore critical points of $I_\eps$ constrained to the fixed point set of the above map, i.e.\ to
\[
  \cO_\eta^{A_N,\tau} = \{(\la,\xi)\in\cO_\eta^{A_N}: \la_i=\dots=\la_1,\ \xi_i=R_k^{i-1}\xi_1, \ i=2,\dots,k\},
\]
are critical points of $I_\eps$. 

In conclusion, for the proofs of Theorems~\ref{thm:k bubbles} and \ref{theorem-ball-k bubbles} it remains to find critical points of $I_\eps$ constrained to $\cO^{A_N,\tau}$ for $\eps>0$ small. The additional symmetry of the solutions stated in Theorem~\ref{theorem-ball-k bubbles}, and also in Theorem~\ref{theorem-ball-k bubbles-add}, is obtained as follows. Since the ball is invariant under the action of $A\in SO(N-2)$ defined by $A(x',x''):=(x',Ax'')$ and since $A*(\la,\xi)=(\la,A\xi)=(\la,\xi)$ for $(\la,\xi)\in \cO_\eta^{A_N}$ it follows from \eqref{eq:sym phi} that
\[
  A*\phi_{\eps,\la,\xi} = \phi_{\eps,\la,A\xi} =\phi_{\eps,\la,\xi}\quad\text{for all $A\in SO(N-2)$,}
\]
hence $u_\eps=V_{\eps,\la,\xi}+\phi_{\eps,\la,\xi}$ satisfies $A*u_\eps=u_\eps$, i.e.\ $u_\eps(x',Ax'')=u_\eps(x',x'')$, for all $A\in SO(N-2)$.

In Theorem~\ref{theorem-ball-k bubbles-add} $V_{\eps,\la,\xi}$ is of the form \eqref{shape of V-add1} with $k=4$. Here $I_\eps$ is invariant under the map
\[
  \wt\tau(\la_1,\la_2,\la_3,\la_4,\la_0,\xi_1,\xi_2,\xi_3,\xi_4) = (\la_3,\la_4,\la_1,\la_2,R_4\xi_4,R_4\xi_1,R_4\xi_2,R_4\xi_3),
\]
so, applying the principle of symmetric criticality once more, a critical point of $I_\eps$ constrained to the set
\[
  \cO_\eta^{A_N,\wt\tau} = \{(\la,\xi)\in\cO_\eta^{A_N}: \la_1=\la_3,\ \la_2=\la_4,\  \xi_i=R_k^{i-1}\xi_1, \ i=2,3,4\}
\]
is an unconstrained critical point of $I_\eps$. This can of course be generalized to any even integer $k\ge4$.

\section{Proof of Theorems~\ref{thm:k bubbles} and \ref{theorem-ball-k bubbles}}\label{sec:proof k=2,3}

In this section we consider $V_{\eps,\la,\xi}=-\sum\limits_{i=1}^{k}PU_{\de_i,\xi_i}+PV_{\mu_\eps,\si_\eps}$ for $k=2$ and $k=3$. The reduced energy is expanded as follows; see \cite[Proposition~5.1]{BarGuo-ANS}.

\begin{lemma}\label{lemma-expansion of J-1}
For $\eps\to0^+$ there holds
\begin{equation*}
I_\eps(\la,\xi) = a_1+a_2\eps-a_3\eps^\al-a_4\eps\ln \eps+\psi(\la,\xi)\eps+o(\eps)
\end{equation*}
$C^1$-uniformly with respect to $(\la,\xi)$ in compact sets of
$\cO_\eta$. The constants are given by $$a_1=\frac{1}{N}(k+1)S_0^{\frac{N}{2}},~
a_2 = \frac{(k+1)}{2^*}\int_{\R^N} U_{1,0}^{2^*}\ln U_{1,0} -
 \frac{k+1}{(2^*)^2}S_0^{\frac{N}2},~
a_3=\frac{1}{2} S_0^{\frac{N-2}{2}}\ov{S}\mu_0, ~a_4=\frac{k+1}{2 \cdot2^*}\int_{\R^N} U_{1,0}^{2^*},$$
where $\ov{S}>0$ is defined by $S_\mu=S_0-\ov{S}\mu +O(\mu^2)$; see \cite[Lemma~A.10]{BarGuo-ANS}. The function $\psi:\cO_\eta\to\R$ is given by
\[
\begin{aligned}
\psi(\la,\xi)
 &= b_1\bigg(H(0,0)\la_0^{N-2} + \sum_{i=1}^kH(\xi_i,\xi_i)\la_i^{N-2}
    + 2\sum_{i=1}^kG(\xi_i,0)\la_i^{\frac{N-2}{2}}\la_0^{\frac{N-2}{2}}\\
 &\hspace{2cm}
   - 2\sum_{i,j=1,i<j}^kG(\xi_i,\xi_j)\la_i^{\frac{N-2}{2}}\la_j^{\frac{N-2}{2}}\bigg)
   - b_2\frac{N-2}{2}\ln(\la_1\la_2\dots\la_{k}\la_0),
\end{aligned}
\]
with 
\[
  b_1=\frac{1}{2}C_0\int_{\R^N}U_{1,0}^{2^*-1},\quad b_2=\frac{1}{2^*}\int_{\R^N} U_{1,0}^{2^*}.
\]
\end{lemma}

It is well known that a stable critical point $(\la,\xi)$ of $\psi$ implies the existence of a critical point $(\la_\eps,\xi_\eps)$ of $I_\eps$ for $\eps>0$ small, and that $(\la_\eps,\xi_\eps)\to(\la,\xi)$ as $\eps\to0$. This applies in particular if $(\la,\xi)$ is an isolated critical point of $\psi$ with nontrivial local degree.

\begin{altproof}{Theorem~\ref{thm:k bubbles}} 
Since the symmetries of $I_\eps$ carry over to $\psi$, for the existence of solutions $u_\eps$ as stated in Theorem~\ref{thm:k bubbles} it is sufficient to find stable critical points of $\psi$ constrained to
\[
  \cO_\eta^{A_N,\tau} = \{(\la,\xi)\in\cO_\eta^{A_N}: \la_i=\la_1,\ \xi_i=R_k^{i-1}\xi_1, \ i=2,\dots,k\},
\]
where $k=2,3$. Observe that $I_\eps$ and $\psi$ are also invariant with respect to the action of $A\in SO(2)$ given by $(x',0)\mapsto (Ax',0)$ acting on the $\xi_i$. Therefore in the case $k=2$, setting $\xi_1=(t,0,\dots,0)$ for $0<t<1$ and $\xi_2=R_2\xi_1=-\xi_1$, it is sufficient to find stable critical points of the function $f_2:\R^+\times\R^+\times(0,1)\to\R$ defined by
\[
\begin{aligned}
f_2(\la_0,\la_1,t)
 &= \psi(\la_0,\la_1,\la_1,\xi_1,-\xi_1)\\
 &= b_1\left(H(0,0)\la_0^{N-2} + 2H(\xi_1,\xi_1)\la_1^{N-2}
        + 4G(\xi_1,0)\la_1^{\frac{N-2}{2}}\la_0^{\frac{N-2}{2}} - 2G(\xi_1,-\xi_1)\la_1^{N-2}\right)\\
 &\hspace{1cm}  - b_2\frac{N-2}{2}\ln\left(\la_1^2\la_0\right).
\end{aligned}
\]
This proves Theorem~\ref{thm:k bubbles} in the case $k=2$. For $k=3$ we set  $\xi_1=(t,0,\dots,0)$ for $0<t<1$, $\xi_2=R_3\xi_1=\left(-\frac{t}{2},\frac{\sqrt{3}t}{2},0,\ldots,0\right)$, $\xi_3=R_3\xi_2=\left(-\frac{t}{2},-\frac{\sqrt{3}t}{2},0,\ldots,0\right)$. As above it is sufficient to find stable critical points of the function $f_3:\R^+\times\R^+\times(0,1)\to\R$ defined by 
\[
\begin{aligned}
f_3(\la_0,\la_1,t)
 &= \psi(\la_0,\la_1,\la_1,\la_1,\xi_1,\xi_2,\xi_3)\\
 &= b_1\left(H(0,0)\la_0^{N-2} + 3H(\xi_1,\xi_1)\la_1^{N-2}
        + 6G(\xi_1,0)\la_0^{\frac{N-2}{2}}\la_1^{\frac{N-2}{2}}
        - 6G(\xi_1,\xi_2)\la_1^{N-2}\right)\\
 &\hspace{1cm}
   - b_2\ln\left(\la_1^3\la_0\right)^{\frac{N-2}{2}}.
\end{aligned}
\]
Here we used that $G(\xi_1,0)=G(\xi_2,0)=G(\xi_3,0)$ and $G(\xi_1,\xi_2)=G(\xi_1,\xi_3)=G(\xi_2,\xi_3)$, as well as $H(\xi_1,\xi_1)=H(\xi_2,\xi_2)=H(\xi_3,\xi_3)$. This proves Theorem~\ref{thm:k bubbles} also in the case $k=3$.
\end{altproof}

\begin{altproof}{Theorem~\ref{theorem-ball-k bubbles}} 
Here we need to find stable critical points of $f_2$ and $f_3$ if $G$ is the Green function of the Dirichlet Laplace operator in the unit ball in $\R^N$. Our proof uses the explicit knowledge of $G$. In the case $k=2$ the proof of \cite[Lemma~3.1]{BarDPis-BLMS} applies almost verbatim and shows that $f_2$ has two isolated critical points: a local saddle point with Morse index $1$, hence with local degree $-1$, and an isolated local minimum, hence with local degree $1$. These are stable critical points, giving rise to critical points of $I_\eps$ for $\eps>0$ small.

In the case $k=3$ we set
\[
\ga_1(t)
 := H(\xi_1,\xi_1)-2G(\xi_1,\xi_2)
   = \frac{1}{(1-t^2)^{N-2}}-\frac{2}{(\sqrt{3}t)^{N-2}}+\frac{2}{(t^4+t^2+1)^{\frac{N-2}{2}}}
\]
and
\begin{equation}\label{eq:tau_1}
\tau_1(t) := G(\xi_1,0) = \frac{1}{t^{N-2}}-1
\end{equation}
so that
\[
  f_3(\la_0,\la_1,t)
    = b_1\left(H(0,0)\la_0^{N-2}+3\ga_1(t)\la_1^{N-2}+6\tau_1(t)\la_0^{\frac{N-2}{2}}\la_1^{\frac{N-2}{2}}\right)
        - b_2\frac{N-2}{2}\ln\left(\la_1^3\la_0\right).
\]
One easily checks that $\ga'_1(t)>0$, $\ga_1(t)\to -\infty$ as $t\to 0^+$, and $\ga_1(\frac12)>0$. Thus there exists $t^*\in(0,\frac{1}{2})$ such that
\begin{equation*}
\ga_1(t^*)=0\quad\text{and}\quad\ga_1(t)>0 \text{ for all } t\in(t^*,1).
\end{equation*}
A direct computation shows that for $t\in(t^*,1)$ there exist unique $\la_0(t)$, $\la_1(t)$ such that
\begin{equation*}
\frac{\pa f_3(\la_0(t),\la_1(t),t)}{\pa\la_0} = 0\quad\text{and}\quad
\frac{\pa f_3(\la_0(t),\la_1(t),t)}{\pa\la_1} = 0.
\end{equation*}
In fact one obtains
\begin{equation}\label{equ-2-add2-lemma blowing up speed}
  \la_0(t)^{\frac{N-2}{2}}=\al(\xi_1,\xi_2)\la_1(t)^{\frac{N-2}{2}}\quad\text{and}\quad
  \la_1(t)^{\frac{N-2}{2}}=\sqrt{\frac{1}{\be(\xi_1,\xi_2)}\cdot\frac{b_2}{2b_1}},
\end{equation}
where
\[
  \al(x,y) = \frac{-2G(x,0)+\sqrt{4G^2(x,0)+4H(0,0)(H(x,x)-2G(x,y))}}{2H(0,0)}
\]
and
\[
  \be(x,y) = H(x,x)-2G(x,y)+G(x,0)\al(x,y).
\]
Moreover, continuing the computation one obtains
\begin{eqnarray*}
\frac{\pa^2 f_3(\la_0(t),\la_1(t),t)}{\pa\la_1^2}
&=& 3(N-2)b_1\left((N-3)\ga_1(t)\la_1^{N-4}
       +\frac{N-4}{2}\tau_1(t)\la_0^{\frac{N-6}{2}}\la_1^{\frac{N-2}{2}}\right)
     +\frac{3(N-2)b_2}{2\la_1^2}\\
&=& 3(N-2)b_1\left((N-2)\ga_1(t)\la_1^{N-4}
         + \frac{N-2}{2}\tau_1(t)\la_0^{\frac{N-6}{2}}\la_1^{\frac{N-2}{2}}\right),\\
\frac{\pa^2 f_3(\la_0(t),\la_1(t),t)}{\pa\la_0^2}
&=& (N-2)b_1\left((N-3)H(0,0)\la_0^{N-4}
        + \frac{3(N-4)}{2}\tau_1(t)\la_0^{\frac{N-2}{2}}\la_1^{\frac{N-6}{2}}\right)\\
&&\hspace{.5cm}  + \frac{(N-2)b_2}{2\la_0^2}\\
&=& (N-2)b_1\left((N-2)H(0,0)\la_0^{N-4}
        +\frac{3(N-2)}{2}\tau_1(t)\la_0^{\frac{N-2}{2}}\la_1^{\frac{N-6}{2}}\right),\\
\frac{\pa^2 f_3(\la_0(t),\la_1(t),t)}{\pa\la_0\pa\la_1}
&=& \frac{3(N-2)^2}{2}b_1\tau_1(t)\la_0^{\frac{N-4}{2}}\la_1^{\frac{N-4}{2}}.
\end{eqnarray*}
It follows that the Hessian matrix $D^2_{\la_0,\la_1}f_3(\la_0(t),\la_1(t),t)$ is positive definite, hence nondegenerate. Therefore it is sufficient to find stable critical points of the function
\[
\nu_1(t)
 := f_3\left(\la_0(t),\la_1(t),t\right)
  = 2b_2-b_2\frac{N-2}{2}\ln\left(\la_1^3(t)\la_0(t)\right).
\]
As in \cite[(3.4)]{BarDPis-BLMS} one sees that
\begin{equation}\label{eq:nu}
 \lim_{t\to (t^*)^+}\nu_1(t)=-\infty\quad\text{and}\quad
\lim_{t\to 1^-}\nu_1(t)=+\infty.
\end{equation}
Now we prove $\nu_1'(\frac{1}{2})<0$ for $N$ large. Set
\[
  \al(t):=\al(\xi_1,\xi_2)=-\tau_1(t)+\sqrt{\tau_1^2(t)+\ga_1(t)},
\]
where we used $H(0,0)=1$. We obtain
\[
\nu_1'(t) = \frac{\pa f_3(\la_0(t),\la_1(t),t)}{\pa t} = 3b_1\big(\ga'_1(t)+2\al(t)\tau'_1(t)\big)\la_1^{N-2}.  
\]
Then setting $\iota_1(t):=\ga'_1(t)+2\al(t)\tau'_1(t)$, it is enough to show $\iota_1\left(\frac12\right)<0$ for $N$ large. In fact, since $\frac{\ga_1(\frac12)}{\tau_1^2(\frac12)}<1$ for $N$ large we see as in \cite[(3.9)]{BarDPis-BLMS} that
\[
\iota_1({\textstyle\frac12})
 \le \ga'_1({\textstyle\frac12})
      + \frac{4\ga_1(\frac12)}{5\tau_1(\frac{1}{2})}\tau'_1({\textstyle\frac12}).
\]
A direct computation gives for $N$ large the inequalites
\[
\ga'_1({\textstyle\frac12})
 = (N-2)\left(({\textstyle\frac43})^{N-1} + 4({\textstyle\frac{2}{\sqrt{3}}})^{N-2}
      - \frac{\frac{3}{2}}{(\frac{1}{16}+\frac{1}{4}+1)^{\frac{N}{2}}}\right)
  < \frac{11(N-2)}{10}({\textstyle\frac43})^{N-1}
\]
and
\[
\tau'_1({\textstyle\frac12}) = -(N-2)2^{N-1}
\]
and
\[
\frac{\ga_1(\frac12)}{\tau_1(\frac12)} = \frac{(\frac{4}{3})^{N-2}-2(\frac{2}{\sqrt{3}})^{N-2}
      + \frac{2}{(\frac{1}{16}+\frac{1}{4}+1)^{\frac{N-2}{2}}}}{2^{N-2}-1}
  > \frac{11}{12}\cdot\frac{(\frac{4}{3})^{N-2}}{2^{N-2}}
\]
which yield $\iota_1(\frac{1}{2})<0$, hence $\nu_1'(\frac{1}{2})<0$ for $N$ large enough. This together with \eqref{eq:nu} implies that $\nu_1$ has a local maximum $t_1\in(t^*,\frac12)$ and a local minimum $t_2\in(\frac12,\infty)$. These are nondegenerate because $\nu_1$ is analytic. In conclusion, $f_3$ has two critical points: $(\la_0(t_1),\la_1(t_1),t_1)$ with Morse index $1$ and $(\la_0(t_2),\la_1(t_2),t_2)$ with Morse index $0$. This concludes the proof of Theorem \ref{theorem-ball-k bubbles}.
\end{altproof}

\begin{remark}\label{rem:nonexist}
a) For $k=3$, $N=7$, numerical computations show that one cannot find $t_0\in(t^*,1)$ such that $\nu_1'(t_0)=0$. Therefore it is necessary to assume $N$ large here.

b) For $k=4$, the idea above cannot give the existence of nodal solutions with five bubbles, one positive at the origin and four negative as in Theorem~\ref{thm:k bubbles}. This is the content of Proposition \ref{prop:k=4}.
\end{remark}

\section{Proof of Proposition~\ref{prop:k=4}}\label{sec:proof prop}

It follows from Lemma~\ref{lemma-expansion of J-1} that $I_\eps$ does not have critical points for $\eps>0$ small if $\psi$ does not have critical points. This also holds if we constrain $I_\eps$ and $\psi$ to $\cO_\eta^{A_N,\tau}$. Setting
$\xi_1=(t,0,\dots,0)$ for $0<t<1$,  $\xi_2=R_4\xi_1=(0,t,0,\dots,0)$, $\xi_3=R_4\xi_2=(-t,0,\dots,0)$, and $\xi_4=R_4\xi_3=(0,-t,\dots,0)$ we need to consider the function
\[
\begin{aligned}
f_4(\la_0,\la_1,t)
 &:= \psi(\la_0,\la_1,\la_1,\la_1,\la_1,\xi_1,\xi_2,\xi_3,\xi_4)\\
 &= b_1\left(H(0,0)\la_0^{N-2} + 4H(\xi_1,\xi_1)\la_1^{N-2} + 8G(\xi_1,0)\la_0^{\frac{N-2}{2}}\la_1^{\frac{N-2}{2}}\right.
\\
&\hspace{1cm}
   -8G(\xi_1,\xi_2)\la_1^{N-2}-4G(\xi_1,\xi_3)\la_1^{N-2} \bigg)-b_2\frac{N-2}{2}\ln(\la_1^4\la_0),
  \end{aligned}
\]
where we use
\[
H(\xi_1,\xi_1) = H(\xi_2,\xi_2) = H(\xi_3,\xi_3)=H(\xi_4,\xi_4) \quad\text{and}\quad
G(\xi_1,0) = G(\xi_2,0) = G(\xi_3,0) = G(\xi_4,0)
\]
as well as
\[
G(\xi_1,\xi_2) = G(\xi_2,\xi_3) = G(\xi_3,\xi_4) = G(\xi_4,\xi_1),\quad
G(\xi_1,\xi_3) = G(\xi_2,\xi_4).
\]
Proposition~\ref{prop:k=4} follows if we can prove that $f_4$ does not have critical points. Let $\tau_1(t)$ be as in \eqref{eq:tau_1} and define
\begin{eqnarray*}
\ga_2(t)
 &:=& H(\xi_1,\xi_1)-2G(\xi_1,\xi_2)-G(\xi_1,\xi_3)\\
 &=& \frac{1}{(1-t^2)^{N-2}} - \frac{2}{(\sqrt{2}t)^{N-2}} + \frac{2}{(t^4+1)^{\frac{N-2}{2}}}
      - \frac{1}{(2t)^{N-2}}+\frac{1}{(t^2+1)^{N-2}}
\end{eqnarray*}
so that
\[
  f_4(\la_0,\la_1,t)
    = b_1\left(H(0,0)\la_0^{N-2}+4\ga_2(t)\la_1^{N-2}+8\tau_1(t)\la_0^{\frac{N-2}{2}}\la_1^{\frac{N-2}{2}}\right)
        - b_2\frac{N-2}{2}\ln\left(\la_1^4\la_0\right).
\]
A direct computation shows that
\[
\ga'_2(t)
 = (N-2)\left(\frac{2t}{(1-t^2)^{N-1}} + \frac{2}{(\sqrt{2})^{N-2}t^{N-1}}
    - \frac{4t^3}{(t^4+1)^{\frac{N}{2}}} + \frac{1}{2^{N-2}t^{N-1}}
    - \frac{2t}{(t^2+1)^{N-1}}\right)
 > 0.
\]
Clearly $\ga_2(t)\to -\infty$ as $t\to 0^+$, and $\ga_2\left(\frac{1}{\sqrt{2}}\right)>0$. Then there exists $t^*\in(0,\frac{1}{\sqrt{2}})$ such that
\begin{equation}\label{equ-1-proof-theorem-ball-k bubbles-2}
\ga_2(t^*)=0 \quad\text{and}\quad \ga_2(t)>0\text{ for all $t\in(t^*,1)$.}
\end{equation}
Notice that 
\begin{equation*}\begin{aligned}
&\la_0^{\frac{N-2}{2}}=\al_1(\xi_1,\xi_2,\xi_3)\la_1^{\frac{N-2}{2}},\quad
\la_1^{\frac{N-2}{2}}=\sqrt{\frac{1}{\be_1(\xi_1,\xi_2,\xi_3)}\cdot\frac{b_2}{2b_1}}, \\
&H(\xi_1,\xi_1)-2G(\xi_1,\xi_2)-G(\xi_1,\xi_3)>0,
\end{aligned}
\end{equation*}
where $$\al_1(x,y,z)=\frac{-3G(x,0)+\sqrt{9G^2(x,0)+4H(0,0)(H(x,x)-2G(x,y)-G(x,z))}}{2H(0,0)},$$
and $$\be_1(x,y,z)=H(x,x)-2G(x,y)-G(x,z)+G(x,0)\al_1(x,y,z).$$
Setting $\al_1(t):=\al_1(\xi_1,\xi_2,\xi_3)=\frac{-3\tau_1(t)+\sqrt{9\tau_1^2(t)+4\ga_2(t)}}{2}$ and $\iota_2(t):=\ga'_2(t)+2\al_1(t)\tau'_1(t)$, a similar argument as above shows  that problem \eqref{pro} admits a solution with $5$ bubbles, one positive at the origin and $4$ negative as in Theorem~\ref{thm:k bubbles} only if $\iota_2(t)$ has a zero in $(t^*,1)$. The following claim implies that this is not the case.

\noindent
{\bf Claim:} If $N\ge7$ then $\iota_2(t)>0$ for any $t\in(t^*,1)$.\\
We first show that  $t^*>\frac{\sqrt{6}-\sqrt{2}}{2}$, where $t^*$ is from  \eqref{equ-1-proof-theorem-ball-k bubbles-2}. In order to see this, it is enough to prove
$\ga_2\left(\frac{\sqrt{6}-\sqrt{2}}{2}\right) < 0$. Since $2^{2/5}\cdot2(\frac{\sqrt{6}-\sqrt{2}}{2})^2<1<(\frac{\sqrt{6}-\sqrt{2}}{2})^4+1$, we have
\[
\frac{1}{(\sqrt{2}\cdot\frac{\sqrt{6}-\sqrt{2}}{2})^{N-2}}
 > \frac{2}{((\frac{\sqrt{6}-\sqrt{2}}{2})^4+1)^{\frac{N-2}{2}}},
\quad\text{for all } N\ge7.
\]
On the other hand, it is easy to see that
\[
\frac{1}{(1-(\frac{\sqrt{6}-\sqrt{2}}{2})^2)^{N-2}}
 = \frac{1}{(\sqrt{2}(\frac{\sqrt{6}-\sqrt{2}}{2}))^{N-2}} \quad\text{and}\quad
\frac{1}{(2(\frac{\sqrt{6}-\sqrt{2}}{2}))^{N-2}}
 > \frac{1}{((\frac{\sqrt{6}-\sqrt{2}}{2})^2+1)^{N-2}}.
\]
It follows that $\ga_2(\frac{\sqrt{6}-\sqrt{2}}{2}) < 0$.

Now we prove $\iota_2(t)>0$ for $t\in(t^*,1)\subset(\frac{\sqrt{6}-\sqrt{2}}{2},1)$. It is easy to see that 
\[
\ga'_2(t) \ge (N-2)\cdot\frac{2t}{(1-t^2)^{N-1}} \quad\text{and}\quad
\ga_2(t) \le \frac{1}{(1-t^2)^{N-2}} \quad\text{for all $t\in\left(\frac{\sqrt{6}-\sqrt{2}}{2},1\right)$.}
\]
Then we have for all $t\in (t^*,1)$ and $N\ge7$:
\begin{eqnarray*}
\frac{\iota_2(t)}{N-2}
 &\ge& \frac{2t}{(1-t^2)^{N-1}}
       - 3\left(\frac{1}{t^{N-2}}-1\right)\cdot\frac{1}{t^{N-1}}
       \left(\sqrt{1+\frac{4\cdot\frac{1}{(1-t^2)^{N-2}}}{9(\frac{1}{t^{N-2}}-1)^2}}-1\right)\\
&\ge& \frac{2t}{(1-t^2)^{N-1}} - \frac{3}{t^{2N-3}}\cdot
      \left(\sqrt{1+\frac{4}{9}(\frac{t^2}{1-t^2})^{N-2}\cdot\frac{1}{(1-t^{N-2})^2}}-1\right).
\end{eqnarray*}
Setting $T:=\frac{t^2}{1-t^2}$ it is enough to prove that
\[
\frac23\cdot T^{N-1}+1
 > \sqrt{1+\frac{4}{9}\cdot T^{N-2}\cdot\frac{1}{(1-t^{N-2})^2}}
\]
which is equivalent to
\begin{equation}\label{equ-3-k=5-impossible-proposition}
(T^N+3T)\cdot(1-t^{N-2})^2>1.
\end{equation}
It is obvious that if $t\in{\textstyle[\frac{1}{\sqrt{2}},\frac45)}$, then
\begin{equation}\label{equ-4-k=5-impossible-proposition}
3T\cdot(1-t^{N-2})^2 \ge 3(1-t^5)^2>1,
\end{equation}
and if $t\in[\frac{4}{5},1)$, then
\begin{equation}\label{equ-5-k=5-impossible-proposition}
T^N\cdot(1-t^{N-2})^2
 \ge T^N\cdot(1-t)^2=\frac{t^4}{(1+t)^2}\cdot T^{N-2}
 > \frac{(\frac{4}{5})^4}{4}(\frac{(\frac{4}{5})^2}{1-(\frac{4}{5})^2})^5
 > 1.
\end{equation}
Now we are left to prove \eqref{equ-3-k=5-impossible-proposition} for $t\in\left(t^*,\frac{1}{\sqrt{2}}\right)$. First of all, if
$t\in\left(\frac{\sqrt{6}-\sqrt{2}}{2},\frac{1}{\sqrt{2}}\right)$, then $T\in\left(\frac{\sqrt{3}-1}{2},1\right)$. Setting
\[
  f(T):=3T\left(1-t^{N-2}\right)^2=3T\left(1-\left(\frac{T}{1+T}\right)^{\frac{N-2}{2}}\right)^2,
\]
a direct computation shows that
\begin{eqnarray*}
f'(T)
 &=& \left(1-\left(\frac{T}{1+T}\right)^{\frac{N-2}{2}}\right)
          \left(3-3\left(\frac{T}{1+T}\right)^{\frac{N-2}{2}}
      - 3(N-2)\left(\frac{T}{1+T}\right)^{\frac{N-2}{2}}\frac{1}{1+T}\right)\\
&\ge& \left(1-\left(\frac12\right)^{\frac{N-2}{2}}\right)
       \left(3-3\left(\frac12\right)^{\frac{N-2}{2}}
       -3(N-2)\left(\frac12\right)^{\frac{N-2}{2}}\frac{1}{1+\frac{\sqrt{3}-1}{2}}\right)\\
&\ge& \left(1-\left(\frac12\right)^{\frac{N-2}{2}}\right)
      \left(3-3\left(\frac12\right)^{\frac52}
       -15\left(\frac12\right)^{\frac52}\frac{1}{1+\frac{\sqrt{3}-1}{2}}\right)
 > 0,
\end{eqnarray*}
where in the second inequality we use the fact that $3-3(\frac{1}{2})^{\frac{N-2}{2}}-3(N-2)(\frac{1}{2})^{\frac{N-2}{2}}\frac{1}{1+\frac{\sqrt{3}-1}{2}}$ is increasing in $N$. Now we conclude that
\begin{eqnarray}\label{equ-6-k=5-impossible-proposition}
f(T)
 > 3\cdot\frac{\sqrt{3}-1}{2}
   \left(1-\left(\frac{\frac{\sqrt{3}-1}{2}}{1+\frac{\sqrt{3}-1}{2}}\right)^{\frac52}\right)^2
 > 1\quad\text{for all } N\ge7.
\end{eqnarray}
The claim, hence Proposition~\ref{prop:k=4}, follows combining \eqref{equ-3-k=5-impossible-proposition}, \eqref{equ-4-k=5-impossible-proposition}, \eqref{equ-5-k=5-impossible-proposition}, and \eqref{equ-6-k=5-impossible-proposition}. 

\section{Proof of Theorem~\ref{theorem-ball-k bubbles-add}}\label{sec:proof 5 bubbles}

In this section we consider solutions of the form $V_{\eps,\la,\xi}=\sum\limits_{i=1}^{k}(-1)^iPU_{\de_i,\xi_i}+PV_\si$. Then the reduced function in Lemma \ref{lemma-expansion
of J-1} becomes
\[
\begin{aligned}
\wt{\psi}(\la,\xi)
 &= b_1\left(H(0,0)\la_0^{N-2}+\sum_{i=1}^kH(\xi_i,\xi_i)\la_i^{N-2}
    + 2\sum_{i=1}^k(-1)^{i-1}G(\xi_i,0)\la_0^{\frac{N-2}{2}}\la_i^{\frac{N-2}{2}}\right.
\\
&\hspace{1cm}
 \left.
  +2\sum_{i,j=1,i<j}^k(-1)^{i+j-1}G(\xi_i,\xi_j)\la_i^{\frac{N-2}2}\la_j^{\frac{N-2}2}\right)
          -b_2\frac{N-2}{2}\ln(\la_0\la_1\la_2\dots\la_{k}),
\end{aligned}
\]
where $b_1, b_2$ are as in Lemma~\ref{lemma-expansion of J-1}.

\begin{altproof}{Theorem~\ref{theorem-ball-k bubbles-add}}
Let $k=4$. Using the symmetry again we set $\xi_1=(t,0,\dots,0)$ for $0<t<1$,  $\xi_2=R_4\xi_1=(0,t,0,\dots,0)$, $\xi_3=R_4\xi_2=(-t,0,\dots,0)$, and $\xi_4=R_4\xi_3=(0,-t,\dots,0)$. As in the proof of Theorem \ref{theorem-ball-k bubbles}, it is sufficient to find stable critical points of $\wt{\psi}$ constrained to
$ \cO_\eta^{A_N,\wt\tau}$. Since 
\[
H(\xi_1,\xi_1)=H(\xi_2,\xi_2)=H(\xi_3,\xi_3)=H(\xi_4,\xi_4),\quad
G(\xi_1,0)=G(\xi_2,0)=G(\xi_3,0)=G(\xi_4,0),
\]
and
\[
G(\xi_1,\xi_2)=G(\xi_2,\xi_3)=G(\xi_3,\xi_4)=G(\xi_4,\xi_1),\quad
G(\xi_1,\xi_3)=G(\xi_2,\xi_4),
\]
we need to find a critical point of the function
\[
\begin{aligned}
f_5(\la_0,\la_1,\la_2,t)
 &:=\wt{\psi}(\la_0,\la_1,\la_2,\la_1,\la_2,\xi_1,\xi_2,\xi_3,\xi_4)\\
 &= b_1\bigg(H(0,0)\la_0^{N-2} + 2H(\xi_1,\xi_1)(\la_1^{N-2}+\la_2^{N-2})
    +4G(\xi_1,0)\la_0^{\frac{N-2}{2}}\left(\la_1^{\frac{N-2}{2}}-\la_2^{\frac{N-2}{2}}\right)\\
&\hspace{1.5cm}
    + 8G(\xi_1,\xi_2)\la_1^{\frac{N-2}{2}}\la_2^{\frac{N-2}{2}}
    - 2G(\xi_1,\xi_3)\la_1^{N-2}-2G(\xi_1,\xi_3)\la_2^{N-2}\bigg)\\
&\hspace{0.5cm}
    - b_2\frac{N-2}{2}\ln\left(\la_0\la_1^2\la_2^2\right).
\end{aligned}
\]

\noindent
{\bf Claim 1:} There exist $t_1^*\in(0,\frac12)$ and $t_2^*\in(\frac12,1)$ such that for $t\in(0,t_1^*)\cup(t_2^*,1)$ the equation 
\begin{equation}\label{eq:nabla f_5}
  \nabla_{\la_0,\la_1,\la_2} f_5(\la_0,\la_1,\la_2,t)=0
\end{equation}
has a unique solution $(\la_0(t),\la_1(t),\la_2(t),t)$.\\
Observe that \eqref{eq:nabla f_5} is equivalent to the equations
\begin{equation}\label{equ-2-ball-k=5-theorem}
H(0,0)\la_0^{N-2} + 2G(\xi_1,0)\la_0^{\frac{N-2}{2}}\left(\la_1^{\frac{N-2}{2}} - \la_2^{\frac{N-2}{2}}\right) = \frac{b_2}{2b_1}
\end{equation}
and
\begin{equation}\label{equ-3-ball-k=5-theorem}
(H(\xi_1,\xi_1) - G(\xi_1,\xi_3))\la_1^{N-2} + G(\xi_1,0)\la_0^{\frac{N-2}{2}}\la_1^{\frac{N-2}{2}}
      + 2G(\xi_1,\xi_2)\la_1^{\frac{N-2}{2}}\la_2^{\frac{N-2}{2}}
  = \frac{b_2}{2b_1},
\end{equation}
and
\begin{equation}\label{equ-4-ball-k=5-theorem}
(H(\xi_1,\xi_1) - G(\xi_1,\xi_3))\la_2^{N-2} - G(\xi_1,0)\la_0^{\frac{N-2}{2}}\la_2^{\frac{N-2}{2}}
       + 2G(\xi_1,\xi_2)\la_1^{\frac{N-2}{2}}\la_2^{\frac{N-2}{2}}
  = \frac{b_2}{2b_1}.
\end{equation}
From \eqref{equ-3-ball-k=5-theorem} and  \eqref{equ-4-ball-k=5-theorem} we deduce
\begin{equation}\label{equ-5-ball-k=5-theorem}
\la_2^{\frac{N-2}{2}}-\la_1^{\frac{N-2}{2}}
 = \frac{G(\xi_1,0)}{H(\xi_1,\xi_1)-G(\xi_1,\xi_3)}\la_0^{\frac{N-2}{2}},
\end{equation}
which combined with \eqref{equ-2-ball-k=5-theorem} implies:
\begin{equation}\label{equ-6-ball-k=5-theorem}
\la_0^{N-2}
 = \frac{H(\xi_1,\xi_1)-G(\xi_1,\xi_3)}
     {H(\xi_1,\xi_1)-G(\xi_1,\xi_3)-2G^2(\xi_1,0)}\cdot\frac{b_2}{2b_1}.
\end{equation}
As a consequence of \eqref{equ-5-ball-k=5-theorem} we get
\[
\la_1^{N-2}+\la_2^{N-2}-2\la_1^{\frac{N-2}{2}}\la_2^{\frac{N-2}{2}}
 = \la_0^{N-2}\cdot\left(\frac{G(\xi_1,0)}{H(\xi_1,\xi_1)-G(\xi_1,\xi_3)}\right)^2
\]
hence using \eqref{equ-3-ball-k=5-theorem} and \eqref{equ-4-ball-k=5-theorem} we deduce:
\begin{equation}\label{equ-7-ball-k=5-theorem}
\la_1^{\frac{N-2}{2}}\la_2^{\frac{N-2}{2}}
 = \frac{1}{H(\xi_1,\xi_1)-G(\xi_1,\xi_3)+2G(\xi_1,\xi_2)}\cdot\frac{b_2}{2b_1}
\end{equation}
and
\begin{equation}\label{equ-8-ball-k=5-theorem}
\la_1^{N-2}+\la_2^{N-2}
 = \frac{1}{H(\xi_1,\xi_1) - G(\xi_1,\xi_3)+2G(\xi_1,\xi_2)}\cdot\frac{b_2}{b_1} + \la_0^{N-2}\cdot\left(\frac{G(\xi_1,0)}{H(\xi_1,\xi_1)-G(\xi_1,\xi_3)}\right)^2.
\end{equation}
Let $\tau_1(t)=G(\xi_1,0)$ be as in \eqref{eq:tau_1} and set
\[
\ga_3(t)
 := H(\xi_1,\xi_1)-G(\xi_1,\xi_3)
  = \frac{1}{(1-t^2)^{N-2}}-\frac{1}{(2t)^{N-2}}+\frac{1}{(t^2+1)^{N-2}}\\
\]
and
\[
\ga_4(t) := G(\xi_1,\xi_2) = \frac{1}{(\sqrt{2}t)^{N-2}}-\frac{1}{(t^4+1)^{\frac{N-2}{2}}}
\]
so that
\[
\begin{aligned}
  f_5(\la_0,\la_1,\la_2,t)
    &= b_1\bigg(H(0,0)\la_0^{N-2}+2\ga_2(t)\left(\la_1^{N-2}+\la_2^{N-2}\right)+8\ga_4(t)\la_1^{\frac{N-2}{2}}\la_2^{\frac{N-2}{2}}\\
    &\hspace{1.5cm}
           +4\tau_1(t)\la_0^{\frac{N-2}{2}}\left(\la_1^{\frac{N-2}{2}}-\la_2^{\frac{N-2}{2}}\right)\bigg)
        - b_2\frac{N-2}{2}\ln\left(\la_1^4\la_0\right).
\end{aligned}
\]

A direct computation shows that $\ga'_3(t)>0$, $\ga_3(t)\to -\infty$ as $t\to 0^+$, $\ga_3(t)\to +\infty$ as $t\to 1^-$, and $\ga_3(\frac{1}{2})>0$. Thus there exists $t_1^*\in(0,\frac{1}{2})$ such that
\[
\ga_3(t_1^*)=0 \quad\text{and}\quad \ga_3(t)<0 \text{ for all } t\in(0,t_1^*).
\]
On the other hand, $\big(\ga_3(t)-2\tau_1^2(t)\big)'>0$, $\ga_3(t)-2\tau_1^2(t)\to -\infty$ as $t\to0^+$, $\ga_3(t)-2\tau_1^2(t)\to +\infty$ as $t\to 1^-$, and $\ga_3(\frac{1}{2})-2\tau_1^2(\frac{1}{2})<0$. Thus there exists $t_2^*\in(\frac{1}{2},1)$ such that
\[
\ga_3(t_2^*)-2\tau_1^2(t_2^*)=0 \quad\text{and}\quad \ga_3(t)-2\tau_1^2(t)>0 \text{ for all } t\in(t_2^*,1).
\]
It follows that for every $t\in(0,t_1^*)\cup(t_2^*,1)$ there exist unique $\la_0(t)$, $\la_1(t)$, $\la_2(t)$ such that
\[
\nabla_{\la_0,\la_1,\la_2} f_5(\la_0(t), \la_1(t),\la_2(t),t)=0,
\]
where $\la_0(t), \la_1(t),\la_2(t)$ satisfy \eqref{equ-5-ball-k=5-theorem}, \eqref{equ-6-ball-k=5-theorem}, \eqref{equ-7-ball-k=5-theorem} and  \eqref{equ-8-ball-k=5-theorem}. This proves Claim 1.

\noindent
{\bf Claim 2:} The Hessian matrix $D^2_{\la_0,\la_1,\la_2}f_5(\la_0(t),\la_1(t),\la_2(t),t)$ is nondegenerate for any $t\in(0,t_1^*)\cup(t_2^*,1)$.\\
A direct computation using  \eqref{equ-2-ball-k=5-theorem}, \eqref{equ-3-ball-k=5-theorem}, and \eqref{equ-4-ball-k=5-theorem} shows that, writing $\la_i$ instead of $\la_i(t)$,
\begin{equation*}
\begin{aligned}
\frac{\pa^2 f_5(\la_0,\la_1,\la_2,t)}{\pa\la_0^2}
 &= (N-2)b_1\bigg((N-3)H(0,0)\la_0^{N-4} + (N-4)\tau_1(t)\la_0^{\frac{N-6}{2}}\left(\la_1^{\frac{N-2}{2}}
      - \la_2^{\frac{N-2}{2}}\right)\bigg)\\
 &\hspace{1cm}
      + \frac{(N-2)b_2}{2\la_0^2}\\
 &= (N-2)^2b_1\bigg(H(0,0)\la_0^{N-4} + \tau_1(t)\la_0^{\frac{N-6}{2}}\left(\la_1^{\frac{N-2}{2}}
      - \la_2^{\frac{N-2}{2}}\right)\bigg),\\
\frac{\pa^2 f_5(\la_0,\la_1,\la_2,t)}{\pa\la_1^2}
 &= (N-2)b_1\bigg(2(N-3)\ga_3(t)\la_1^{N-4}
     +(N-4)\tau_1(t)\la_0^{\frac{N-2}{2}}\la_1^{\frac{N-6}{2}}\\
 &\hspace{1cm}
      + 2(N-4)\ga_4(t)\la_1^{\frac{N-6}{2}}\la_2^{\frac{N-2}{2}}\bigg)
       +\frac{(N-2)b_2}{\la_1^2}\\
 &= (N-2)^2b_1\bigg(2\ga_3(t)\la_1^{N-4}
     + \tau_1(t)\la_0^{\frac{N-2}{2}}\la_1^{\frac{N-6}{2}}
     + 2\ga_4(t)\la_1^{\frac{N-6}{2}}\la_2^{\frac{N-2}{2}}\bigg),\\
\frac{\pa^2 f_5(\la_0,\la_1,\la_2,t)}{\pa\la_2^2}
 &= (N-2)b_1\bigg(2(N-3)\ga_3(t)\la_2^{N-4}
     - (N-4)\tau_1(t)\la_0^{\frac{N-2}{2}}\la_2^{\frac{N-6}{2}}\\
 &\hspace{1cm}
      + 2(N-4)\ga_4(t)\la_1^{\frac{N-6}{2}}\la_2^{\frac{N-2}{2}}\bigg)
       + \frac{(N-2)b_2}{\la_2^2}\\
 &= (N-2)^2b_1\bigg(2\ga_3(t)\la_2^{N-4}
     - \tau_1(t)\la_0^{\frac{N-2}{2}}\la_2^{\frac{N-6}{2}}
     + 2\ga_4(t)\la_1^{\frac{N-2}{2}}\la_2^{\frac{N-6}{2}}\bigg),\\
\frac{\pa^2 f_4(\la_0,\la_1,\la_2,t)}{\pa\la_0\pa\la_1}
  &= (N-2)^2b_1\tau_1(t)\la_0^{\frac{N-4}{2}}\la_1^{\frac{N-4}{2}},\\
\frac{\pa^2 f_4(\la_1,\la_2,\la_0,t)}{\pa\la_0\pa\la_2}
  &=-(N-2)^2b_1\tau_1(t)\la_0^{\frac{N-4}{2}}\la_2^{\frac{N-4}{2}},\\
\frac{\pa^2 f_4(\la_1,\la_2,\la_0,t)}{\pa\la_1\pa\la_2}
  &= 2(N-2)^2b_1\ga_4(t)\la_1^{\frac{N-4}{2}}\la_2^{\frac{N-4}{2}}.
\end{aligned}
\end{equation*}
For simplicity, we introduce the notation
\[
X:=\la_0^{\frac{N-2}{2}},\quad Y:=\la_1^{\frac{N-2}{2}},\quad Z:=\la_2^{\frac{N-2}{2}}.
\]
In order to prove that $D^2_{\la_0,\la_1,\la_2}f_5(\la_0(t),\la_1(t),\la_2(t),t)$ is nondegenerate for any $t\in(0,t_1^*)\cup(t_2^*,1)$, it suffices to show that the matrix
\[
\left(
\begin{array}{ccc}
X+\tau_1(t)(Y-Z) & \tau_1(t)X^\frac{2}{N-2}Y^{\frac{N-4}{N-2}} & -\tau_1(t)X^\frac{2}{N-2}Z^{\frac{N-4}{N-2}}\\
\tau_1(t)X^{\frac{N-4}{N-2}}Y^\frac{2}{N-2}  & 2\ga_3(t)Y+\tau_1(t)X+2\ga_4(t)Z & 2\ga_4(t)Y^{\frac{2}{N-2}}Z^\frac{N-4}{N-2}\\
-\tau_1(t)X^{\frac{N-4}{N-2}}Z^\frac{2}{N-2} & 2\ga_4(t)Y^\frac{N-4}{N-2}Z^{\frac{2}{N-2}} & 2\ga_3(t)Z-\tau_1(t)X+2\ga_4(t)Y \\
\end{array}
\right)
\]
is nondegenerate. Using \eqref{equ-2-ball-k=5-theorem}, \eqref{equ-3-ball-k=5-theorem} and \eqref{equ-4-ball-k=5-theorem} this is equivalent to showing that the matrix
\[
\left(
\begin{array}{ccc}
\frac{X}{2}+\frac{b_2}{4b_1}\cdot\frac{1}{X}& \tau_1(t)X^\frac{2}{N-2}Y^{\frac{N-4}{N-2}} & -\tau_1(t)X^\frac{2}{N-2}Z^{\frac{N-4}{N-2}}\\
\tau_1(t)X^{\frac{N-4}{N-2}}Y^\frac{2}{N-2}  & \ga_3(t)Y+\frac{b_2}{2b_1}\cdot\frac{1}{Y} & 2\ga_4(t)Y^{\frac{2}{N-2}}Z^\frac{N-4}{N-2}\\
-\tau_1(t)X^{\frac{N-4}{N-2}}Z^\frac{2}{N-2} & 2\ga_4(t)Y^\frac{N-4}{N-2}Z^{\frac{2}{N-2}} & \ga_3(t)Z+\frac{b_2}{2b_1}\cdot\frac{1}{Z} \\
\end{array}
\right)
\]
is nondegenerate. A direct computation, using \eqref{equ-7-ball-k=5-theorem}, shows that the determinant of the above matrix has the same sign as $\ga_3(t)$, hence is nontrivial, proving Claim 2.

Theorem~\ref{theorem-ball-k bubbles-add} now follows from

\noindent
{\bf Claim 3:} The function $\nu_2(t):=f_5\big(\la_0(t),\la_1(t),\la_2(t),t\big)$ has a critical point $t_1\in(0,t_1^*)$.\\
Observe that, writing again $\la_i$ instead of $\la_i(t)$,
\[
\begin{aligned}
\nu_2'(t) &= \frac{\pa f_4(\la_0(t),\la_1(t),\la_2(t),t)}{\pa t}\\
 &= 2b_1\bigg(\ga'_3(t)\left(\la_1^{N-2} + \la_2^{N-2}\right)
    + 2\tau'_1(t)\la_0^{\frac{N-2}2}\left(\la_1^{\frac{N-2}2}-\la_2^{\frac{N-2}2}\right)
    + 4\ga'_4(t)\la_1^{\frac{N-2}{2}}\la_2^{\frac{N-2}{2}}\bigg),
\end{aligned}
\]
where $\la_0,\la_1,\la_2$ satisfy \eqref{equ-5-ball-k=5-theorem}, \eqref{equ-6-ball-k=5-theorem}, \eqref{equ-7-ball-k=5-theorem} and \eqref{equ-8-ball-k=5-theorem}. Therefore, $\nu_2'(t)=0$ for $t\in(0,t_1^*)$ is equivalent to
\[
\begin{aligned}
\iota_3(t)
 &:= \ga'_3(t)\big(2\ga_3(t)(\ga_3(t) - 2\tau_1^2(t)) + \tau_1^2(t)(\ga_3(t)+2\ga_4(t))\big)
       - 2\tau'_1(t)\tau_1(t)\ga_3(t)\big(\ga_3(t)+2\ga_4(t)\big)\\
 &\hspace{1cm}
       + 4\ga'_4(t)\ga_3(t)\big(\ga_3(t)-2\tau_1^2(t)\big)\\
 &= 0.
\end{aligned}
\]
It is easy to check that $\iota_3(t)\to-\infty$ as $t\to 0^+$ and $\iota_3(t_1^*)>0$ because $\ga'_3(t_1^*)>0$, $\ga_4(t_1^*)>0$ and $\ga_3(t_1^*)=0$. Hence there exists $t_1\in(0,t_1^*)$ such that $\iota_3(t_1)=0$. Claim 3 follows, finishing the proof of Theorem~\ref{theorem-ball-k bubbles-add}.
\end{altproof}

\begin{remark}
We conjecture that there should also exist $t_2\in (t_2^*,1)$ such that $\iota_3(t_2)=0$. This is not considered here because the computations get enormous.
\end{remark}

\noindent\textbf{Acknowledgements:} The authors would like to thank Professor Daomin Cao for many helpful discussions during the preparation of this paper. This work was carried out while Qianqiao Guo was visiting Justus-Liebig-Universit\"{a}t Gie{\ss}en, to which he would like to
express his gratitude for their warm hospitality.

\noindent\textbf{Funding:} Qianqiao Guo was supported  by the National Natural Science Foundation of China (Grant No. 11971385) and the Natural Science Basic Research Plan in Shaanxi Province of China (Grant No. 2019JM275).


E-mail:\\ 

thomas.bartsch@math.uni-giessen.de

gqianqiao@nwpu.edu.cn

\end{document}